\documentclass[12pt,a4paper]{amsart}
\usepackage{geometry}
\usepackage{tabls}
\usepackage{url}
\usepackage[utf8]{inputenc}
\usepackage{amsfonts}
\usepackage{amssymb}
\usepackage{mathrsfs}
\usepackage{amsmath}
\usepackage{amsthm}
\usepackage{amscd}
\usepackage{fancyhdr}
\usepackage{enumerate}
\usepackage{paralist}
\usepackage{graphicx}
\usepackage{cite}
\usepackage{rotating}

\usepackage{footmisc}
\numberwithin{equation}{section}

\theoremstyle{plain}

\newtheorem{Thm}[equation]{Theorem}
\newtheorem{lem}[equation]{Lemma}
\newtheorem{prop}[equation]{Proposition}
\newtheorem{rem}[equation]{Remark}
\newtheorem{ex}[equation]{Example}

\begin{document}

\title{A continuous version of multiple zeta functions and multiple zeta values}

\author{Jiangtao Li}

\email{lijiangtao@csu.edu.cn}
\address{Jiangtao Li \\ School of Mathematics and Statistics, HNP-LAMA, Central South University, Hunan Province, China}

\begin{abstract}
In this paper we  define a continuous version of multiple zeta functions. They can be analytically continued to  meromorphic functions on $\mathbb{C}^r$ with only simple poles at some special hyperplanes.  The evaluations of these functions at positive integers (continuous multiple zeta values) satisfy the shuffle product. We give a detailed analysis about the depth structure of continuous multiple zeta values. There are also  sum formulas for continuous multiple zeta values. Lastly we calculate some special continuous multiple zeta values in terms of special values of multiple polylogarithms.
\end{abstract}

\let\thefootnote\relax\footnotetext{
2020 $\mathnormal{Mathematics} \;\mathnormal{Subject}\;\mathnormal{Classification}$. 11M32.\\
$\mathnormal{Keywords:}$  Riemann zeta function, Multiple zeta values. }

\maketitle

\section{Introduction}
    For $r\geq 1$, the multiple zeta function is defined by 
     $$\zeta(s_1, s_2,\cdots s_r)=\sum_{0<n_1<n_2<\cdots<n_r}\frac{1}{n_1^{s_1}n_2^{s_2}\cdots n_r^{s_r}}.$$  
     For $(s_1,s_2.\cdots, s_r)=(k_1,k_2,\cdots,k_r)$, $k_1,\cdots, k_{r-1}\geq 1, k_r\geq 2$, the values 
     \[
     \zeta(k_1,k_2,\cdots,k_r)
     \]
     are called multiple zeta values. Multiple zeta values satisfy the stuffle product and the  shuffle product \cite{ikz}. As a result, there are many relations among multiple zeta values.
     
     For $r=1, s_1=1$, it is well-known that the harmonic series
     \[
     \sum_{n\geq 1} \frac{1}{n}
     \]
     is divergent. But the modified version     
      \[
     \sum_{1\leq n\leq k}\frac{1}{n} -\int^{k+1}_1\frac{dx}{x}
     \]
     is convergent to a real number which is called Euler constant as $k\to +\infty.$
     The number $\sum\limits_{1\leq k\leq n}\frac{1}{n}$ is a discrete sum, while the  integral $\int^{k+1}_1\frac{dx}{x}$ can be viewed as a continuous sum. Thus, in some sense, the Euler constant can be viewed as a difference value between a discrete sum and a continuous sum.
     
     From the following formula
     \[
     \begin{split}
     &\;\;\;\;\frac{1}{n_1^{s_1}(n_1+n_2)^{s_2}\cdots (n_1+\cdots+n_r)^{s_r}}\\
     &=\mathop{\int\cdots\int}_{[n_1,n_1+1]\times [n_2,n_2+1]\times \cdots [n_r,n_r+1]}\frac{dx_1dx_2\cdots dx_r}{[x_1]^{s_1}([x_1]+[x_2])^{s_2}\cdots ([x_1]+\cdots+[x_r])^{s_r}},
     \end{split}
     \]
 the multiple zeta functions can be viewed as 
     \[
     \begin{split}
     &\;\;\;\; \zeta(s_1,s_2,\cdots, s_r)\\
     &=\sum_{n_1,\cdots,n_r\geq 1}\frac{1}{n_1^{s_1}(n_1+n_2)^{s_2}\cdots (n_1+\cdots+n_r)^{s_r}}\\
     &= \sum_{n_1,\cdots,n_r\geq 1}  \mathop{\int\cdots\int}_{[n_1,n_1+1]\times [n_2,n_2+1]\times \cdots [n_r,n_r+1]}\frac{dx_1dx_2\cdots dx_r}{[x_1]^{s_1}([x_1]+[x_2])^{s_2}\cdots ([x_1]+\cdots+[x_r])^{s_r}}\\ &=\mathop{\int\cdots\int}_{[1,+\infty)^r}\frac{dx_1dx_2\cdots dx_r}{[x_1]^{s_1}([x_1]+[x_2])^{s_2}\cdots ([x_1]+\cdots+[x_r])^{s_r}},\\
     \end{split}
     \]
     where $[x]$ denotes the Gauss rounding function.
     
     Inspired by the above observations, we define the continuous version of multiple zeta functions as
     \[
     \zeta^{\mathcal{C}}(s_1,s_2,\cdots, s_r)=\mathop{\int\cdots\int}_{[1,+\infty)^r}\frac{dx_1dx_2\cdots dx_r}{x_1^{s_1}(x_1+x_2)^{s_2}\cdots (x_1+\cdots+x_r)^{s_r}}.
     \]
     The multiple zeta functions are constructed from the discrete function $[x]$.  The continuous multiple zeta functions are constructed from the continuous function $x$.      
          We have 
     \begin{Thm}\label{cmzv}
     The continuous multiple zeta function  $ \zeta^{\mathcal{C}}(s_1,s_2,\cdots, s_r)$     is convergent for 
     \[
     \mathrm{Re}(s_1+s_2+\cdots+s_r)>r,\; \mathrm{Re}(s_2+\cdots+s_r)>r-1,\;\cdots, \mathrm{Re}(s_r)>1.
     \]
     Moreover, it can be analytically continued to a meromorphic function on $\mathbb{C}^r$ with possible poles at some special hyperplanes.      \end{Thm}
    The detailed structure of the possible poles of $ \zeta^{\mathcal{C}}(s_1,s_2,\cdots, s_r)$ will be given in Section \ref{acc}.

      The classical multiple zeta values satisfy the stuffle product and the shuffle product \cite{ikz}. For the continuous multiple zeta values (values of the continuous multiple zeta functions at positive integers), one has 
      \begin{Thm}\label{shu} (Rough version)
      The continuous multiple zeta values satisfy the shuffle product.
      \end{Thm} 
      From the definition of the Euler constant,  we also expect that by studying the algebra of multiple zeta values and the algebra of continuous multiple zeta values together, one can find an appropriate algebra to understand the mysterious Euler constant.
      
   Denote by $\mathcal{Z}^{\mathcal{C}}$ the $\mathbb{Q}$-linear space generated by $1$ and  continuous multiple zeta values.   
        By Theorem \ref{shu}, the $\mathbb{Q}$-linear vector space $\mathcal{Z}^{\mathcal{C}}$   is actually a $\mathbb{Q}$-algebra. Thus it is interesting to investigate the structure of this algebra.   
   
    For the continuous multiple zeta values $\zeta^{\mathcal{C}}(k_1,\cdots, k_r)$, the number $r$ is called  its depth.
     For $r=1,k\geq 1$, it is easy to see that
      \[
      \zeta^{\mathcal{C}}(1+k)=\int_1^{+\infty}\frac{dx}{x^{1+k}}=\frac{1}{k}.
      \]
          In general cases, we have 
     \begin{Thm}\label{dep}
     For $m_1,\cdots,m_s\geq 1$, define
            \[
            \zeta^{\mathcal{C}}_{m_1,m_2,\cdots,m_s}(\underbrace{1,\cdots,1,2}_s)=\int\limits^{+\infty}_{m_1}\cdots \int\limits_{m_s}^{+\infty}\frac{dx_1 \cdots dx_{s-1} dx_s}{x_1\cdots(x_1+\cdots +x_{s-1}) (x_1+\cdots+x_s)^2}.
            \]
     Denote by $\mathfrak{D}_r\mathcal{Z}^{\mathcal{C}}$ the $\mathbb{Q}$-linear space generated by the continuous multiple zeta values of depth $ r$, then one has \\
      $ (i)$  $ \mathfrak{D}_1\mathcal{Z}^{\mathcal{C}}\subseteq \mathfrak{D}_2\mathcal{Z}^{\mathcal{C}} \subseteq\cdots \subseteq \mathfrak{D}_r\mathcal{Z}^{\mathcal{C}}\subseteq \cdots;    $\\
   $(ii)$ $  \mathfrak{D}_r\mathcal{Z}^{\mathcal{C}}\subseteq \langle \zeta^{\mathcal{C}}_{m_1,\cdots,m_s}(\underbrace{1,\cdots,1,2}_{s})\,|\,s\geq 1, m_1+\cdots+m_s=r, m_1,\cdots,m_s\geq 1      \rangle_{\mathbb{Q}}$,
     where the right side denotes the $\mathbb{Q}$-linear space generated by the following elements
     \[
     \zeta^{\mathcal{C}}_{m_1,\cdots,m_s}(\underbrace{1,\cdots,1,2}_{s}),\,s\geq 1, m_1+\cdots+m_s=r, m_1,\cdots,m_s\geq 1;
          \]
          $(iii)$ As a result, $$\mathrm{dim}_{\mathbb{Q}} \mathfrak{D}_r\mathcal{Z}^{\mathcal{C}}\leq 2^{r-1}, \forall\,r\geq 1.$$     
     \end{Thm}
     
     For the depth structure of classical multiple zeta values, there is the well-known Broadhurst-Kreimer conjecture \cite{bk}, which is related to the generating series of cusp forms of $\mathrm{SL}_2(\mathbb{Z})$.  By Theorem \ref{dep}, the depth structure of continuous multiple zeta values is much simpler.
     
     For the classical multiple zeta values, one has the sum formulas: 
     \[
     \sum_{\substack{k_1+\cdots+k_{r-1}+k_r=k\\ k_1,\cdots,k_{r}\geq 1}}\zeta(k_1,\cdots, k_{r-1},1+k_r)=\zeta(1+k).
     \]
     For the continuous multiple zeta values, we have
     \begin{Thm}\label{sum}
       For $r\geq 2, k>2(r-1)$, denote by 
     \[
    f(x_1,x_2,\cdots,x_r)=x_r(x_{r-1}+x_{r}-2)\cdots [x_2+\cdots +x_r-2(r-2)]  [x_1+\cdots +x_r-2(r-1)].
     \]
     Denote by ${\bf V}$ the $\mathbb{Q}$-linear space generated by 
     \[
     \frac{1}{x^l},\frac{1}{(x+1)^l},\cdots, \frac{1}{(x+n)^l},\cdots, \forall\, l\geq 1.
     \]
     Define $\eta: {\bf V}\rightarrow {\bf V}$ as the $\mathbb{Q}$-linear transformation which satisfies 
     \[
     \eta\left( \frac{1}{(x+n)^l}    \right)=\frac{1}{n+1}\left( \frac{1}{x^l}-\frac{1}{(x+n+1)^l}\right), \forall \,n\geq 0, l\geq 1.
     \]
     Then
     \[
     \begin{split}
     &\sum_{\substack{k_1+\cdots+k_r=k\\ k_1,\cdots, k_r\geq 1}}f(k_1,\cdots,k_{r-1},k_r)   \zeta^{\mathcal{C}}(k_1,\cdots, k_{r-1},1+k_r)
     =  \underbrace{\eta\circ\cdots \circ \eta}_{r-1}\left( \frac{1}{x^{l}}   \right)  \Bigg{|}_{x=1}, \\
     \end{split}
     \]
     where $l=k-2(r-1)$ and  $\alpha(x)\big{|}_{x=t}$ means $\alpha(t)$ for any rational function $\alpha(x)$.
     
     \end{Thm}
     
     The theory of multiple polylogarithms is related to the arithmetic theory of number fields.
     At the end of this paper, we discuss the relation between continuous multiple zeta values and evaluations of multiple polylogarithms. We also discuss the depth defect phenomena in the algebra of continuous multiple zeta values.  The results in Section 3 and Section 4 reveal that there are interesting  relations between continuous multiple zeta values and cyclotomic multiple zeta values.

      \section{Analytic continuation of continuous multiple zeta functions}\label{acc}
      In this section we show that the continuous multiple zeta functions are convergent under some natural conditions. Furthermore, they can be analytically continued to meromorphic functions with only simple poles at some special hyperplanes.       
      \begin{prop}\label{conv} For $r\geq 1$, the continuous multiple zeta function 
       \[
     \zeta^{\mathcal{C}}(s_1,s_2,\cdots, s_r)=\mathop{\int\cdots\int}_{[1,+\infty)^r}\frac{dx_1dx_2\cdots dx_r}{x_1^{s_1}(x_1+x_2)^{s_2}\cdots (x_1+\cdots+x_r)^{s_r}}
     \]
       is convergent if
           \[
     \mathrm{Re}(s_1+s_2+\cdots+s_r)>r,\; \mathrm{Re}(s_2+\cdots+s_r)>r-1,\;\cdots, \mathrm{Re}(s_r)>1.
     \]
      \end{prop}
         \noindent{\bf Proof:} Denote by $s_l=\sigma_l+it_l, \sigma_l,t_l\in \mathbb{R}, 1\leq l\leq r$. If 
         \[
         \sigma_1+\sigma_2+\cdots+\sigma_r>r, \sigma_2+\cdots+\sigma_r>r-1,\cdots, \sigma_r>1,
         \]
         for $M>1$, we have
         \[
         \begin{split}
         &\;\;\;\; \Bigg{|} \mathop{\int\cdots\int}_{[1,M]^r}\frac{dx_1dx_2\cdots dx_r}{x_1^{s_1}(x_1+x_2)^{s_2}\cdots (x_1+\cdots+x_r)^{s_r}}\Bigg{|}
       \\
         &< \mathop{\int\cdots\int}_{[1,M]^r}  \Bigg{|}     \frac{1}{x_1^{s_1}(x_1+x_2)^{s_2}\cdots (x_1+\cdots+x_r)^{s_r}} \Bigg{|}    dx_1dx_2\cdots dx_r       \\
         &= \mathop{\int\cdots\int}_{[1,M]^r}   \frac{dx_1dx_2\cdots dx_r}{x_1^{\sigma_1}(x_1+x_2)^{\sigma_2}\cdots (x_1+\cdots+x_r)^{\sigma_r}} \\
         &= \frac{1}{\sigma_r-1} \mathop{\int\cdots\int}_{[1,M]^{r-1}}  \frac{dx_1dx_2\cdots dx_{r-1}}{x_1^{\sigma_1}(x_1+x_2)^{\sigma_2}\cdots (x_1+\cdots+x_{r-1})^{\sigma_{r-1}} (1+x_1+\cdots+x_{r-1})^{\sigma_r-1}}             \\
         &<  \frac{1}{\sigma_r-1}   \mathop{\int\cdots\int}_{[1,M]^{r-1}}  \frac{dx_1dx_2\cdots dx_{r-1}}{x_1^{\sigma_1}(x_1+x_2)^{\sigma_2}\cdots (x_1+\cdots+x_{r-1})^{\sigma_{r-1}+\sigma_r-1} }              \\
         &<\frac{1}{(\sigma_r-1)(\sigma_r+\sigma_{r-1}-2)\cdots (\sigma_r+\cdots+\sigma_1-r)}.         \\  \end{split}
         \]
         Since the above inequality holds for any $M>1$, the limit
         $$\lim_{M\rightarrow +\infty} \mathop{\int\cdots\int}_{[1,M]^r}   \frac{dx_1dx_2\cdots dx_r}{x_1^{\sigma_1}(x_1+x_2)^{\sigma_2}\cdots (x_1+\cdots+x_r)^{\sigma_r}} $$
         is convergent.  As a result, the continuous  multiple zeta function $$ \zeta^{\mathcal{C}}(s_1,s_2,\cdots, s_r)$$ is convergent. $\hfill\Box$\\  
         
         The following lemma will be useful in the analytic continuation of continuous multiple zeta function.
         \begin{lem}\label{an}
         If $\varphi(t)$ is an infinitely differentiable bounded function on $(-\epsilon, 1+\epsilon)$ for some $\epsilon>0$, then
         \[
         I_{\varphi}(s)=\int_0^1\varphi(t)t^{s-1}dt
        \]
        can be analytically continued to a meromorphic function on $\mathbb{C}$ with only simple poles at $s=0,-1,\cdots,-n,\cdots$ and $\mathrm{Res}_{s=-n}I_{\varphi}(s)=\frac{ \varphi^{(n)}(0)}{n!}  $.
         \end{lem}
          \noindent{\bf Proof:}   For $k\geq 1$, $\mathrm{Re}(s)>0$, we have
          \[
          I_{\varphi}(s)=\int^1_0\left(\varphi(t)-\sum_{n=0}^k\frac{\varphi^{(n)}(0)}{n!}t^n   \right)t^{s-1}dt+\sum_{n=0}^k\frac{\varphi^{(n)}(0)}{n!}\frac{1}{s+n}.   \]    
          Denote by $R_{\varphi}(t)=\varphi(t)-\sum_{n=0}^k\frac{\varphi^{(n)}(0)}{n!}t^n $. As 
          \[
          R_{\varphi}(t)=\mathrm{O}(t^k),\,t\to 0,
          \]
          it follows that
          \[
          \int^1_0             R_{\varphi}(t)      t^{s-1}dt\]
           is a holomorphic function for $\mathrm{Re}(s)>-k$. 
           
           Since the above analysis holds for any $k\geq 1$. The lemma is proved.  $\hfill\Box$\\        
                      
         For the continuous multiple zeta function  $\zeta^{\mathcal{C}}(s_1,s_2,\cdots, s_r)$, we have 
         \[
         \begin{split}
         &\;\;\;\;\;\zeta^{\mathcal{C}}(s_1,s_2,\cdots, s_r)\\
         &=\mathop{\int\cdots\int}_{[0,1]^r}\frac{1}{\left(\frac{1}{x_1}\right)^{s_1} \left(\frac{1}{x_1}+\frac{1}{x_2}\right)^{s_2} \cdots \left(\frac{1}{x_1}+\cdots+\frac{1}{x_r}\right)^{s_r}   }  \frac{dx_1}{ x_1^2 }\frac{dx_2}{x_2^2 }\cdots \frac{dx_r}{ x_r^2}    \\
         &=\sum_{\sigma\in S_r} \mathop{\int\cdots \int}_{0<x_{\sigma(1)}<x_{\sigma(2)}<\cdots<x_{\sigma(r)}<1} \frac{1}{\left(\frac{1}{x_1}\right)^{s_1} \left(\frac{1}{x_1}+\frac{1}{x_2}\right)^{s_2} \cdots \left(\frac{1}{x_1}+\cdots+\frac{1}{x_r}\right)^{s_r}   }   \frac{dx_1}{ x_1^2 }\frac{dx_2}{x_2^2 }\cdots \frac{dx_r}{ x_r^2}   \\
         &=  \sum_{\sigma\in S_r} \mathop{\int\cdots \int}_{0<x_{1}<x_{2}<\cdots<x_{r}<1} \frac{1}{\left(\frac{1}{x_{\sigma(1)}}\right)^{s_1} \left(\frac{1}{x_{\sigma(1)}}+\frac{1}{x_{\sigma(2)}}\right)^{s_2} \cdots \left(\frac{1}{x_{\sigma(1)}}+\cdots+\frac{1}{x_{\sigma(r)}}\right)^{s_r}  }\frac{dx_1}{ x_1^2 }\frac{dx_2}{x_2^2 }\cdots \frac{dx_r}{ x_r^2},     
         \end{split}
         \]    
         where $S_r$ is the permutation group of the set $\{1,2, \cdots,r\}$. For $\sigma\in S_r$, denote by
         \[
         \begin{split}
         &\;\;\;\; I_{\sigma}(s_1,s_2,\cdots, s_r)\\
         &= \mathop{\int\cdots \int}_{0<x_{1}<x_{2}<\cdots<x_{r}<1} \frac{1}{\left(\frac{1}{x_{\sigma(1)}}\right)^{s_1} \left(\frac{1}{x_{\sigma(1)}}+\frac{1}{x_{\sigma(2)}}\right)^{s_2} \cdots \left(\frac{1}{x_{\sigma(1)}}+\cdots+\frac{1}{x_{\sigma(r)}}\right)^{s_r}  }\frac{dx_1}{ x_1^2 }\frac{dx_2}{x_2^2 }\cdots \frac{dx_r}{ x_r^2}, \\
                  \end{split}
         \]
         one has
         \[
         \zeta^{\mathcal{C}}(s_1,s_2,\cdots, s_r)=\sum_{\sigma\in S_r}     I_{\sigma}(s_1,s_2,\cdots, s_r).               \]
         
         \begin{prop}\label{ac}
         For a fixed $\sigma\in S_r$, define $$m_i=\mathop{\mathrm{min}}_{1\leq j\leq i}\{\sigma(j)\}, \;\forall\, 1\leq i\leq r.$$ 
         Then  $I_{\sigma}(s_1,s_2,\cdots,s_r)$ can be analytically continued to a meromorphic function on $\mathbb{C}^r$ with possible poles at 
         \[
           m_1s_1=2-k_1, m_1s_1+m_2s_2=3-k_2,\cdots,m_1s_1+\cdots+m_rs_r=(r+1)-k_r,\]
           where $ k_1,k_2,\cdots,k_r\geq 1$.   
                    \end{prop}
          \noindent{\bf Proof:}  It is clear that
          \[
          r\geq m_1\geq m_2\geq \cdots\geq m_r=1.
          \] 
          We have
          \[
          \begin{split}
          &\;\;\;\; I_{\sigma}(s_1,s_2,\cdots, s_r)\\
                    &=  \mathop{\int\cdots \int}_{0<x_{1}<x_{2}<\cdots<x_{r}<1} x_1^{m_1s_1-2}x_2^{m_2s_2-2}\cdots x_r^{m_rs_r-2} \varphi(x_1,x_2,\cdots,x_r)dx_1dx_2\cdots dx_r,       \\
          \end{split}
          \]
          where 
          \[
          \varphi(x_1,x_2,\cdots,x_r)= \frac{1}{\left(\frac{x_{m_1}}{x_{\sigma(1)}}\right)^{s_1} \left(\frac{x_{m_2}}{x_{\sigma(1)}}+\frac{x_{m_2}}{x_{\sigma(2)}}\right)^{s_2} \cdots \left(\frac{x_{m_r}}{x_{\sigma(1)}}+\cdots+\frac{x_{m_r}}{x_{\sigma(r)}}\right)^{s_r}  }.
                              \]
                              By letting \[
                              x_1=y_1y_2\cdots y_r, x_2=y_2\cdots y_r, x_r=y_r,
                              \]
                              we have
          \[
          \begin{split}
          &\;\;\;\; I_{\sigma}(s_1,s_2,\cdots, s_r)\\
                    &=  \mathop{\int\cdots \int}_{[0,1]^r} y_1^{m_1s_1-2}y_2^{m_ss_1+m_2s_2-3}\cdots y_r^{m_1s_1+\cdots+m_rs_r-(r+1)} \Phi(y_1,y_2,\cdots,y_r)dy_1dy_2\cdots dy_r  ,       \\
          \end{split}
          \]
          where 
          \[
          \Phi(y_1,y_2,\cdots,y_r)=\varphi(x_1x_2\cdots x_r, x_2\cdots x_r, \cdots, x_r) .
                              \]
           By the definition of $(m_1,m_2,\cdots,m_r)$, it follows that $\Phi(y_1,y_2,\cdots,y_r)$ is a infinitely differentiable function on $(-\epsilon,1+\epsilon)^r$ for some $\epsilon>0$ and $\Phi(0,0,\cdots,0)=1$.
           
           By Lemma \ref{an},  it follows that $I_{\sigma}(s_1,s_2,\cdots,s_r)$ can be analytically continued to a meromorphic function on $\mathbb{C}^r$ with possible poles at the following hyperplanes:
           \[
           m_1s_1=2-k_1, m_1s_1+m_2s_2=3-k_2,\cdots,m_1s_1+\cdots+m_rs_r=(r+1)-k_r,\]
           where $ k_1,k_2,\cdots,k_r\geq 1$.          $\hfill\Box$\\       
           
           By Proposition \ref{conv} and Proposition \ref{ac}, Theorem \ref{cmzv} is proved. 
           
           \begin{rem}
           Zhao \cite{zhao} proved that the multiple zeta function $$\zeta(s_1, s_2, \cdots, s_r)$$ can be analytically continued to a meromorphic function on $\mathbb{C}^r$ with possible poles at some special hyperplanes. Proposition \ref{ac} shows that the structure of the poles of the continuous multiple zeta function $$ \zeta^{\mathcal{C}}(s_1,s_2,\cdots, s_r)$$ is more complicated than that of the multiple zeta function $$\zeta(s_1, s_2, \cdots, s_r).$$           \end{rem}
                           
      \section{Continuous multiple zeta values}
      In this section, firstly we will prove that the continuous multiple zeta values satisfy the shuffle product. Here we will use the notations in \cite{ikz}.
      Secondly, we will give a detailed analysis of the depth structure of continuous multiple zeta values. Lastly, we will show that there are also sum formulas for continuous multiple zeta values.
      
      \subsection{The algebra of continuous multiple zeta values}
           Define $\mathfrak{H}=\mathbb{Q}\langle x,y\rangle$ as the non-commutative polynomial ring over $\mathbb{Q}$ in two indeterminates $x$ and $y$. The shuffle product $\rotatebox{90}{$\rotatebox{180}{$\exists$}$}$ on $\mathfrak{H}$ is defined by 
           \[1\,\rotatebox{90}{$\rotatebox{180}{$\exists$}$}  \,         w=w\, \rotatebox{90}{$\rotatebox{180}{$\exists$}$}\,1=w,\]
           \[uw_1\, \rotatebox{90}{$\rotatebox{180}{$\exists$}$}   \, vw_2=u( w_1\, \rotatebox{90}{$\rotatebox{180}{$\exists$}$}   \, vw_2)+ v( uw_1\, \rotatebox{90}{$\rotatebox{180}{$\exists$}$}   \, w_2),    \]
      for any $u,v \in \{x,y\}$ and $w, w_1,w_2\in \mathfrak{H}$ inductively. Under the shuffle product $\rotatebox{90}{$\rotatebox{180}{$\exists$}$}$, $\mathfrak{H}$ is a commutative  $\mathbb{Q}$-algebra.  Denote by $\mathfrak{H}_{ \tiny{   \rotatebox{90}{$\rotatebox{180}{$\exists$}$}    }  }$ this commutative $\mathbb{Q}$-algebra. Let $\mathfrak{H}^0=\mathbb{Q}+y\mathfrak{H}x$. It is clear that $\mathfrak{H}^0$ is a $\mathbb{Q}$-subalgebra of $\mathfrak{H}_{ \tiny{   \rotatebox{90}{$\rotatebox{180}{$\exists$}$}    }  }$.
      
      Define a $\mathbb{Q}$-linear map by $$Z: \mathfrak{H}^0\rightarrow \mathcal{Z}^{\mathcal{C}},$$
      \[
      Z(1)=1, Z(yx^{k_1-1}yx^{k_2-1}\cdots yx^{k_r-1})=\zeta^{\mathcal{C}}(k_1,k_1,\cdots,k_r),
      \]
      for $k_1,\cdots,k_{r-1},\geq 1, k_r\geq 2$. The precise version of Theorem \ref{shu} is 
 \begin{Thm}\label{shuf}
 The $\mathbb{Q}$-linear map $Z:\mathfrak{H}^0\rightarrow \mathcal{Z}^{\mathcal{C}}$ is an algebra homomorphism under the shuffle product $\rotatebox{90}{$\rotatebox{180}{$\exists$}$}$  on $\mathfrak{H}^0$.
 \end{Thm}
  \noindent{\bf Proof:}  For convenience, let $k_0=0$.
  For $k_1,\cdots,k_{r-1},\geq 1, k_r\geq 2$, it is easy to check that
  \[   
  \frac{1}{x_1^{k_1}(x_1+x_2)^{k_2}\cdots (x_1+\cdots+x_r)^{k_r}}=\mathop{\int\cdots\int}_{+\infty >t_1>t_2>\cdots>t_k>0}\omega_1(t_1)\omega_2(t_2)\cdots \omega_k(t_k).
  \]
  Here $k=k_1+k_2+\cdots+k_r$ and
   \[
   \omega_{k_0+\cdots+k_j+1}(t_{k_0+\cdots+k_j+1})=e^{-x_{j+1} t_{k_0+\cdots+k_j+1}} dt_{k_0+\cdots+k_j+1}, \;0\leq j\leq r-1, \]
 \[
 \omega_l(t_l)=dt_l,\; l\neq k_0+\cdots +k_j+1, \forall\, \;0\leq j\leq r-1.  \]
  Thus 
  \[
  \begin{split}
  &\;\;\;\;\zeta^{\mathcal{C}}(k_1,k_2,\cdots,k_r)\\
  &=\mathop{\int\cdots\int}_{[1,+\infty)^r}\frac{dx_1dx_2\cdots dx_r}{x_1^{k_1}(x_1+x_2)^{k_2}\cdots (x_1+\cdots+x_r)^{k_r}}\\
  &=\mathop{\int\cdots\int}_{[1,+\infty)^r}\left(  \mathop{\int\cdots\int}_{+\infty>t_1>t_2>\cdots>t_k>0}\omega_1(t_1)\omega_2(t_2)\cdots \omega_k(t_k)     \right)dx_1dx_2\cdots dx_r\\
  &= \mathop{\int\cdots\int}_{+\infty>t_1>t_2>\cdots>t_k>0}\Omega_1(t_1)\Omega_2(t_2)\cdots \Omega_k(t_k),   \end{split}
  \]
  where 
  \[
  \begin{split}
  &\;\;\;\;\Omega_{k_0+\cdots+k_j+1}( t_{k_0+\cdots+k_j+1})\\
  &=\int_1^{+\infty}\left(e^{-x_{j+1} t_{k_0+\cdots+k_j+1}} dt_{k_0+\cdots+k_j+1}\right)dx_{j+1}=\frac{e^{- t_{k_0+\cdots+k_j+1}}}{t_{k_0+\cdots+k_j+1}} dt_{k_0+\cdots+k_j+1}, 0\leq j\leq r-1,\\
    \end{split}
  \]
  and 
  \[
   \Omega_l(t_l)=dt_l,\; l\neq k_0+\cdots +k_j+1, \forall\, 1\leq l\leq k.  \]
   By the theory of iterated path integrals \cite{chen}, the theorem is proved.   $\hfill\Box$\\    
   
    \begin{rem}
 By letting $u_1=e^{-t_1},u_2=e^{-t_2},\cdots, u_k=e^{-t_k}$ in the above theorem, one has 
 \[
 \zeta^{\mathcal{C}}(k_1,k_2,\cdots,k_r)=\mathop{\int\cdots\int}_{0<u_1<\cdots<u_k<1}\lambda_1(u_1)\lambda_2(u_2)\cdots \lambda_k(u_k),
 \]
 where 
\[
\lambda_i(u)=
\begin{cases}
&\frac{du}{\mathrm{ln}\frac{1}{u}},\;i\in \{1, k_1+1,\cdots, k_1+\cdots+k_{r-1}+1\},\\
&\frac{du}{u},\; i\notin \{1, k_1+1,\cdots, k_1+\cdots+k_{r-1}+1\}.\\
\end{cases}
\]
On the other hand, the classical multiple zeta values are defined by
 \[
 \zeta(k_1,k_2,\cdots,k_r)=\mathop{\int\cdots\int}_{0<t_1<\cdots<t_k<1}\omega_1(t_1)\omega_2(t_2)\cdots \omega_k(t_k),
 \]
 where 
\[
\omega_i(t)=
\begin{cases}
&\frac{dt}{1-t},\;i\in \{1, k_1+1,\cdots, k_1+\cdots+k_{r-1}+1\},\\
&\frac{dt}{t},\; i\notin \{1, k_1+1,\cdots, k_1+\cdots+k_{r-1}+1\}.\\
\end{cases}
\]
 \end{rem}

      \subsection{The depth structure of continuous multiple zeta values}\label{de}
      The depth structure of multiple zeta values is extremely complicated. For the depth two and depth three cases, now there are only the expected upper bound of the dimension of depth-graded version multiple zeta values (see \cite{bk}, \cite{gkz} and \cite{gon}). For the higher depth cases, it is still a myth.
      
     In this subsection, we will give a detailed analysis of the depth structure of continuous multiple zeta values. For convenience, for any set $\mathcal{A}\subseteq \mathbb{R}$, denote by $\langle \mathcal{A}\rangle_{\mathbb{Q}}$ the $\mathbb{Q}$-linear space generated by the elements in $\mathcal{A}$.
    
    By definition,    for $r\geq 1$,   $\mathfrak{D}_r\mathcal{Z}^{\mathcal{C}}$ is the $\mathbb{Q}$-linear space generated by the continuous multiple zeta values of depth $ r$. We wish to prove that     $\mathfrak{D}_r\mathcal{Z}^{\mathcal{C}} \subseteq \mathfrak{D}_{r+1}\mathcal{Z}^{\mathcal{C}}$. This statement follows immediately from the following observation:
    \[
    \begin{split}
    &\;\;\;\;\zeta^{\mathcal{C}}(k_1,k_2,\cdots,k_r)\\
    &=\mathop{\int\cdots\int}_{[1,+\infty)^r}\frac{dx_1dx_2\cdots dx_r}{x_1^{k_1}(x_1+x_2)^{k_2}\cdots (x_1+\cdots+x_r)^{k_r}}\\
    &=\mathop{\int\cdots\int}_{[1,+\infty)^{r}}\frac{dx_1dx_2\cdots dx_{r}}{x_1^{k_1}(x_1+x_2)^{k_2}\cdots (x_1+\cdots+x_{r})^{k_{r}-1}  ( x_1+\cdots+x_{r}     +1  )  }    \\
    &\;\;\;\;\;\;\;\;\;\;\;\;\;\;\;\;\;+\mathop{\int\cdots\int}_{[1,+\infty)^{r}}\frac{dx_1dx_2\cdots dx_{r}}{x_1^{k_1}(x_1+x_2)^{k_2}\cdots (x_1+\cdots+x_{r})^{k_{r}}  ( x_1+\cdots+x_{r}     +1  )  }  \\
    &=\zeta^{\mathcal{C}}(k_1,k_2,\cdots,k_{r-1},k_r-1,2)+\zeta^{\mathcal{C}}(k_1,k_2,\cdots,k_{r-1},k_r,2).
                    \end{split}
     \]
     
     For $r=1$, the statement $(ii)$ of Theorem \ref{dep} is obvious. Now we assume that $r\geq 2$.
     For $\zeta^{\mathcal{C}}(k_1,k_2,\cdots,k_r)$, one has
     \[
     \zeta^{\mathcal{C}}(k_1,k_2,\cdots,k_r)=\int^{+\infty}_1\frac{F(x_1)}{x_1^{k_1}}dx_1,
     \]
     where $$ F(x_1)=\mathop{\int\cdots\int}\limits_{[1,+\infty)^{r-1}} \frac{dx_2\cdots dx_r}{(x_1+x_2)^{k_2}\cdots (x_1+\cdots+x_r)^{k_r}} .$$
      For $k_1\geq 2$, by integration by parts, we have
      \[
        \zeta^{\mathcal{C}}(k_1,k_2,\cdots,k_r)=\frac{F(1)}{k_1-1}+\frac{1}{k_1-1}\int_1^{+\infty}\frac{F^\prime(x_1)}{x_1^{k_1-1}}dx_1.
              \]
              By induction, it follows that 
              \[
              \begin{split}
              &\;\;\;\; \zeta^{\mathcal{C}}(k_1,k_2,\cdots,k_r)\\
              &=\frac{F(1)}{k_1-1}+\frac{F^\prime(1)}{(k_1-1)(k_1-2)}+\cdots+\frac{F^{(k-2)}(1)}{(k_1-1)!}+\frac{1}{(k_1-1)!}\int_1^{+\infty}\frac{F^{(k-2)}(x_1)}{x_1}dx_1.
              \end{split}
                             \]
                  For $m_1,m_2,\cdots, m_r\geq 1$, define
                  \[
                  \zeta^{\mathcal{C}}_{m_1,m_2,\cdots,m_r}(k_1,k_2,\cdots, k_r)=\mathop{\int\cdots\int}_{H_{m_1,m_2,\cdots,m_r}}\frac{dx_1dx_2\cdots dx_r}{x_1^{k_1}(x_1+x_2)^{k_2}\cdots (x_1+\cdots+x_r)^{k_r}},
                  \] 
                  where $H_{m_1,m_2,\cdots,m_r}=[m_1,+\infty)\times [m_2,+\infty)\times \cdots \times [m_r,+\infty)$.  
                  
                  From the above analysis, it is clear that
               \[
               \zeta^{\mathcal{C}}(k_1,k_2,\cdots,k_r)\in \langle  \zeta^{\mathcal{C}}_{2,1,\cdots,1} (l_2,\cdots, l_r), \zeta^{\mathcal{C}}(1,l_2,\cdots,l_r)\,|\,l_2,\cdots,l_{r-1}\geq 1, l_r\geq 2\,  \rangle_{\mathbb{Q}}.
               \]  
               For $l_2\geq 2$, one can use the similar trick to
               \[
               \zeta^{\mathcal{C}}_{2,1,\cdots,1} (l_2,\cdots, l_r), \zeta^{\mathcal{C}}(1,l_2,\cdots,l_r)\,,l_2,\cdots,l_{r-1}\geq 1, l_r\geq 2.                              \]
                Thus we have   
               \[
              \begin{split}
               &\zeta^{\mathcal{C}}(k_1,k_2,\cdots,k_r)\in
                \langle  \zeta^{\mathcal{C}}_{3,1,\cdots,1} (l_3,\cdots, l_r), \zeta^{\mathcal{C}}_{2,1,\cdots,1} (1,l_3,\cdots, l_r), \\
                &\;\;\;\;\;\;\;\;\;\;\;\;\;\;\;\;\;\;\;\;\;\;\; \zeta^{\mathcal{C}}_{1,2,\cdots,1}(1,l_3,\cdots,l_r),\zeta^{\mathcal{C}}(1,1,l_3,\cdots,l_r)
               \,|\,l_3,\cdots,l_{r-1}\geq 1, l_r\geq 2\,  \rangle_{\mathbb{Q}}.\\
               \end{split}
               \]  
               By repeating the above procedure, at last we have 
               $$  \mathfrak{D}_r\mathcal{Z}^{\mathcal{C}}\subseteq \langle \zeta^{\mathcal{C}}_{m_1,\cdots,m_s}(\underbrace{1,\cdots,1,2}_{s})\,|\,s\geq 1, m_1+\cdots+m_s=r, m_1,\cdots,m_s\geq 1      \rangle_{\mathbb{Q}}.$$
               As a result, the statement $(ii)$ of Theorem \ref{dep} is proved.
               
               The statement $(iii)$ of Theorem \ref{dep} follows from 
               \[
               \begin{split}
               &\;\;\;\;\mathrm{dim}_{\mathbb{Q}}  \mathfrak{D}_r\mathcal{Z}^{\mathcal{C}}\\
               &\leq \sum_{s=1}^r \mathrm{dim}_{\mathbb{Q}} \langle \zeta^{\mathcal{C}}_{m_1,\cdots,m_s}(\underbrace{1,\cdots,1,2}_{s})\,|\,m_1+\cdots+m_s=r, m_1,\cdots,m_s\geq 1      \rangle_{\mathbb{Q}}   \\
               &\leq  \sum_{s=1}^r  \sharp \big{\{} (m_1,\cdots,m_s)\,\big{|}\, m_1+\cdots+m_s=r, m_1,\cdots, m_s\geq 1                   \big{\}}\\
               &=\sum_{s=1}^r\binom{r-1}{s-1}\\
               &=2^{r-1},
               \end{split}
               \]
               for $r\geq 1$. Here $\sharp \mathcal{A}$ means the number of elements of $\mathcal{A}$ for any finite set $\mathcal{A}$.

      Now we give some explicit calculations about continuous multiple zeta values.
      \begin{lem} For $x_1,\cdots, x_r\neq 0$, one has
      \[
  \frac{1}{x_1x_2\cdots x_r}=\sum_{\sigma\in S_r} \frac{1}{x_{\sigma(1)}(x_{\sigma(1)}+x_{\sigma(2)})\cdots (x_{\sigma(1)}+x_{\sigma(2)}+\cdots x_{\sigma(r)})},
    \]
         \end{lem}
            \noindent{\bf Proof:}  For $1\leq i\leq r$, one can check that
            \[
            \frac{1}{x_i}=\int^{+\infty}_1e^{-x_it_i}dt_i.
            \]
            Thus 
            \[
            \begin{split}
            &\;\;\;\; \frac{1}{x_1x_2\cdots x_r}\\
            &=\mathop{\int\cdots\int}_{[1,+\infty)^r}   e^{-(x_1t_1+x_2t_2+\cdots+x_rt_r)}dt_1\cdots dt_r        \\
            &=\sum_{\sigma\in S_r}\mathop{\int\cdots\int}_{+\infty>t_{\sigma(1)}>t_{\sigma(2)}>\cdots>t_{\sigma(r)}>1}     e^{-(x_1t_1+x_2t_2+\cdots+x_rt_r)}dt_1\cdots dt_r   \\
            &=\sum_{\sigma\in S_r}\mathop{\int\cdots\int}_{+\infty>t_{1}>t_{2}>\cdots>t_{r}>1}     e^{-(x_1t_{\sigma(1)}+x_2t_{\sigma(2)}+\cdots+x_rt_{\sigma(r)})}dt_1\cdots dt_r              \\
            &=\sum_{\sigma\in S_r}\mathop{\int\cdots\int}_{+\infty>t_{1}>t_{2}>\cdots>t_{r}>1}  e^{-(x_{\sigma(1)}t_{1}+x_{\sigma(2)}t_{2}+\cdots+x_{\sigma(r)}t_{r})}dt_1\cdots dt_r  \\
            &=\sum_{\sigma\in S_r} \frac{1}{x_{\sigma(1)}} \mathop{\int\cdots\int}_{+\infty>t_{2}>\cdots>t_{r}>1}  e^{-[(x_{\sigma(1)}+x_{\sigma(2)})t_{2}+x_{\sigma(3)} t_3+\cdots+x_{\sigma(r)}t_{r}]}dt_1\cdots dt_r  \\
            &\;\;\;\;\;\;\;\;\;\cdots\;\;\;\;\;\cdots \;\;\;\;\;\cdots\\
            &=\sum_{\sigma\in S_r} \frac{1}{x_{\sigma(1)}(x_{\sigma(1)}+x_{\sigma(2)})\cdots (x_{\sigma(1)}+x_{\sigma(2)}+\cdots x_{\sigma(r)})}.
            \end{split}
            \]
            $\hfill\Box$\\               
                  \begin{prop}\label{algebra}
      (i) For $r=2,k_1,k_2\geq 2$, we have 
    \[\zeta^{\mathcal{C}}(1,2)=\mathrm{log}\,2,\zeta^{\mathcal{C}}(k_1,k_2)\in\mathbb{Q}\,\mathrm{log}\,2+\mathbb{Q} \]
      where $\mathbb{Q}\, \mathrm{log}\,2+\mathbb{Q} $ denotes the $\mathbb{Q}$-linear space generated by $\mathrm{log}\,2, 1$;\\
      (ii) For $r=3,k_1,k_2\geq 1,k_3\geq 2$, we have 
           \[
      \zeta^{\mathcal{C}}(k_1,k_2,k_3)\in  \mathbb{Q}\,\zeta^{\mathcal{C}}(1,1,2)   +  \mathbb{Q}\,\mathrm{log}\,3+   \mathbb{Q}\,\mathrm{log}\,2+\mathbb{Q},    \]
            where $$ \mathbb{Q}\,\zeta^{\mathcal{C}}(1,1,2)   +  \mathbb{Q}\,\mathrm{log}\,3+  \mathbb{Q}\, \mathrm{log}\,2+\mathbb{Q} $$ denotes the $\mathbb{Q}$-linear space generated by $ \zeta^{\mathcal{C}}(1,1,2), \mathrm{log}\,3,    \mathrm{log}\,2, 1$;   \\
            (iii) For $r\geq 2$,
            \[
            \zeta^{\mathcal{C}}(\underbrace{1,\cdots,1}_{r-1},2)= \mathop{\int\cdots \int}_{[0,1]^{r-1}} \frac{dy_1 dy_2\cdots dy_{r-1}}{1+y_1+\cdots+y_1\cdots y_{r-1}}  .            \]                          \end{prop}
  \noindent{\bf Proof:}          
  The statements $(i)$
  and $(ii)$       follow immediately from Theorem \ref{dep}.
  Since 
  \[
  \frac{1}{x_1x_2\cdots x_r}=\sum_{\sigma\in S_r} \frac{1}{x_{\sigma(1)}(x_{\sigma(1)}+x_{\sigma(2)})\cdots (x_{\sigma(1)}+x_{\sigma(2)}+\cdots x_{\sigma(r)})},
    \]
    we have
    \[
    \begin{split}
    &\;\;\;\;\mathop{\int\cdots \int}_{[1,+\infty)^r}\frac{dx_1dx_2\cdots dx_r}{x_1x_2\cdots x_r(x_1+x_2+\cdots+x_r)}\\
    &=\sum_{\sigma\in S_r} \mathop{\int\cdots \int}_{[1,+\infty)^r}    \frac{dx_1dx_2\cdots dx_r }{x_{\sigma(1)}(x_{\sigma(1)}+x_{\sigma(2)})\cdots (x_{\sigma(1)}+x_{\sigma(2)}+\cdots x_{\sigma(r)})^2}       \\
    &=r!\zeta^{\mathcal{C}}(\underbrace{1,\cdots,1}_{r-1},2).
        \end{split}
    \]
    Thus 
    \[
    \begin{split}
    &\;\;\;\;\zeta^{\mathcal{C}}(\underbrace{1,\cdots,1}_{r-1},2 )\\
    &=\frac{1}{r!} \mathop{\int\cdots \int}_{[1,+\infty)^r}\frac{dx_1dx_2\cdots dx_r}{x_1x_2\cdots x_r(x_1+x_2+\cdots+x_r)}\\
    &=  \frac{1}{r!} \mathop{\int\cdots \int}_{[0,1]^r}\frac{dx_1dx_2\cdots dx_r}{x_1x_2\cdots x_r(\frac{1}{x_1}+\frac{1}{x_2}+\cdots+\frac{1}{x_r})}\\
    &= \mathop{\int\cdots \int}_{0<x_1<x_2<\cdots<x_r<1}\frac{dx_1dx_2\cdots dx_r}{x_1x_2\cdots x_r(\frac{1}{x_1}+\frac{1}{x_2}+\cdots+\frac{1}{x_r})}.\\
    \end{split}
    \]
    By letting $x_1=y_1y_2\cdots y_r, x_2=y_2\cdots y_r, \cdots, x_r=y_r$, one has 
    \[
    \begin{split}
    &\;\;\;\;\zeta^{\mathcal{C}}(\underbrace{1,\cdots,1}_{r-1},2 )\\
    &=\mathop{\int\cdots \int}_{[0,1]^r}\frac{y_2\cdots y_r^{r-1}dy_1dy_2\cdots dy_r}{y_1y_2^2\cdots y_r^r  (\frac{1}{y_1y_2\cdots y_r}+\frac{1}{y_2\cdots y_r}+\cdots +\frac{1}{y_r}               )                     }\\
    &=\mathop{\int\cdots \int}_{[0,1]^r} \frac{dy_1 dy_2\cdots dy_r}{1+y_1+\cdots+y_1y_2\cdots y_{r-1}}\\
    &=\mathop{\int\cdots \int}_{[0,1]^{r-1}} \frac{dy_1 dy_2\cdots dy_{r-1}}{1+y_1+\cdots+y_1y_2\cdots y_{r-1}}.\\
                \end{split}
    \]
$\hfill\Box$\\   
            
\begin{rem}
For $r=1$, it is clear that $$\mathrm{dim}_{\mathbb{Q}}  \mathfrak{D}_1\mathcal{Z}^{\mathcal{C}}=1.$$ For $r=2$,
by the Hermite-Lindemann Transcendence Theorem, $\mathrm{log}\,2$ is a transcendental number. Thus $$\mathrm{dim}_{\mathbb{Q}}  \mathfrak{D}_2\mathcal{Z}^{\mathcal{C}}=2.$$

Unfortunately, in most cases (for example $r\geq 4$), we only have $$\mathrm{dim}_{\mathbb{Q}}  \mathfrak{D}_r\mathcal{Z}^{\mathcal{C}}<2^{r-1}.$$
The above inequality follows from the following simple observation:
\[
\langle \zeta^{\mathcal{C}}_{m_1,m_2}(1,2)\,|\, m_1+m_2=r, m_1,m_2\geq 1\rangle_{\mathbb{Q}}=\langle  \mathrm{log}\,r,\mathrm{log}\,(r-1),\cdots, \mathrm{log}\,2  \rangle_{\mathbb{Q}}.
\]
For $r\geq 4$, it is clear that 
\[
\mathrm{dim}_{\mathbb{Q}}\langle  \mathrm{log}\,r,\mathrm{log}\,(r-1),\cdots, \mathrm{log}\,2  \rangle_{\mathbb{Q}}<r-1.
\]
     \end{rem}  
     
     \begin{rem}
     For $r=1,2,3$, one can check that 
     $$  \mathfrak{D}_r\mathcal{Z}^{\mathcal{C}}= \langle \zeta^{\mathcal{C}}_{m_1,\cdots,m_s}(\underbrace{1,\cdots,1,2}_{s})\,|\,s\geq 1, m_1+\cdots+m_s=r, m_1,\cdots,m_s\geq 1      \rangle_{\mathbb{Q}}.$$
     It is not clear that whether this statement is true or not for $r\geq 4$.
     \end{rem}                      
        
        \begin{rem}
        From Yamamoto \cite{yam}, it follows that
 \[
 \mathop{\int}_{[0,1]^{2k}}\frac{dx_1dx_2\cdots dx_{2k}}{1-x_1+\cdots +(-1)^ix_1\cdots x_i+\cdots+x_1\cdots x_{2k}}=\zeta^{\star}(\underbrace{2,\cdots,2}_{k})
, \]
\[
\mathop{\int}_{[0,1]^{2k+1}}\frac{dx_1dx_2\cdots dx_{2k+1}}{1-x_1+\cdots +(-1)^ix_1\cdots x_i+\cdots-x_1\cdots x_{2k+1}}=\zeta^{\star}(1,\underbrace{2,\cdots,2}_{k})
.\]
Here $\zeta^{\star}(k_1,k_2,\cdots,k_r)$ denotes the multiple zeta-star values
\[
\zeta^{\star}(k_1,k_2,\cdots,k_r)=\sum_{0<n_1\leq n_2\leq \cdot\leq n_r}\frac{1}{n_1^{k_1}n_2^{k_2}\cdots n_r^{k_r}}.
\]
One can compare the above formulas with Proposition \ref{algebra}, (iii).
        \end{rem}

   \subsection{Sum formulas for continuous multiple zeta values} Continuous multiple zeta values also satisfy some kinds of sum formulas. The essential reason is due to the following lemma.
   \begin{lem}\label{slem}
   For $K\geq 2$, one has 
   \[
   \sum_{\substack{n_1+n_2=K\\n_1,n_2\geq 1}}\frac{1}{x^{n_1}(x+c)^{n_2}}=\frac{1}{c}\left[\frac{1}{x^{K-1}}-\frac{1}{(x+c)^{K-1}}\right].
  \]
   \end{lem}
    \noindent{\bf Proof:}  We have
    \[
    \begin{split}
    &\;\;\;\;\sum_{\substack{n_1+n_2=K\\n_1,n_2\geq 1}}\frac{1}{x^{n_1}(x+c)^{n_2}}\\
    &= \sum_{\substack{n_1+n_2=K\\n_1,n_2\geq 1}}\frac{1}{c} \frac{(x+c)-x}{x^{n_1}(x+c)^{n_2}}\\
    &=\frac{1}{c} \sum_{\substack{n_1+n_2=K\\n_1,n_2\geq 1}} \left(\frac{1}{x^{n_1}(x+c)^{n_2-1}} -\frac{1}{x^{n_1-1}(x+c)^{n_2}}    \right) \\
    &=\frac{1}{c}\left(\frac{1}{x^{K-1}}+ \sum_{\substack{n_1+n_2=K-1\\n_1,n_2\geq 1}}\frac{1}{x^{n_1}(x+c)^{n_2}}   \right)-\frac{1}{c}\left( \sum_{\substack{n_1+n_2=K-1\\n_1,n_2\geq 1}}\frac{1}{x^{n_1}(x+c)^{n_2}}+ \frac{1}{(x+c)^{K-1}}  \right)    \\
    &=\frac{1}{c} \left[\frac{1}{x^{K-1}}-\frac{1}{(x+c)^{K-1}}\right] .  \end{split}
        \]    $\hfill\Box$\\

    Now we are ready to prove Theorem \ref{sum}. We have
       \[
     \begin{split}
     &\;\;\;\;\sum_{\substack{k_1+\cdots+k_r=k\\ k_1,\cdots, k_r\geq 1}}f_r(k_1,\cdots,k_{r-1},k_r)   \zeta^{\mathcal{C}}(k_1,\cdots, k_{r-1},1+k_r)\\
    & =[k-2(r-1)] \cdot \\
    &\sum_{\substack{k_1+\cdots+k_r=k\\ k_1,\cdots, k_r\geq 1}}k_r(k_{r-1}+k_{r}-2)\cdots [k_2+\cdots +k_r-2(r-2)]  \zeta^{\mathcal{C}}(k_1,\cdots, k_{r-1},1+k_r)  \\
    &=[k-2(r-1)]  
    \sum_{\substack{k_1+\cdots+k_r=k\\ k_1,\cdots, k_r\geq 1}}(k_{r-1}+k_{r}-2)\cdots [k_2+\cdots +k_r-2(r-2)]   \cdot \\
    &\;\;\;\;\;\;\;\;\;\;\;\;\;\;\;\;\; \mathop{\int \cdots\int}_{[1,+\infty)^{r-1}} \frac{dx_1\cdots dx_{r-1}}{x_1^{k_1}\cdots (x_1+\cdots+x_{r-1})^{k_{r-1}}  (x_1+\cdots+x_{r-1}+1)^{k_{r}}   }  . 
     \end{split}
     \]
     By Lemma \ref{slem}, for $k-2(r-1)>0$, one has 
     \[
     \begin{split}
     &\;\;\;\;\sum_{\substack{k_1+\cdots+k_r=k\\ k_1,\cdots, k_r\geq 1}}f(k_1,\cdots,k_{r-1},k_r)   \zeta^{\mathcal{C}}(k_1,\cdots, k_{r-1},1+k_r)\\
     \end{split}
     \]
     \[
     \begin{split}
    & =[k-2(r-1)] 
    \sum_{\substack{k_1+\cdots+k_{r-2}+l_{r-1}=k\\k_1,\cdots,k_{r-2}\geq 1,l_{r-1}\geq 2}}(l_{r-1}-2)\cdots (k_2+\cdots+k_{r-2}+l_{r-1}-2(r-2))\cdot\\
    &\mathop{\int\cdots\int}_{[1,+\infty)^{r-1}}
    \frac{1}{x_1^{k_1}\cdots (x_1+\cdots +x_{r-2})^{k_{r-2}}}\left(\frac{dx_1\cdots dx_{r-2}dx_{r-1} }{(x_1+\cdots+x_{r-1})^{l_{r-1}-1}}-   \frac{dx_1\cdots dx_{r-2}dx_{r-1} }{(x_1+\cdots+x_{r-1}+1)^{l_{r-1}-1}}   \right)\\
      &  =[k-2(r-1)] 
    \sum_{\substack{k_1+\cdots+k_{r-2}+l_{r-1}=k\\k_1,\cdots,k_{r-2}\geq 1,l_{r-1}> 2}}(l_{r-1}-2)\cdots (k_2+\cdots+k_{r-2}+l_{r-1}-2(r-2))\cdot\\
 &\mathop{\int\cdots\int}_{[1,+\infty)^{r-1}}
    \frac{1}{x_1^{k_1}\cdots (x_1+\cdots +x_{r-2})^{k_{r-2}}}\left(\eta\left(\frac{1}{x^{l_{r-1}-1}}\right)\Bigg{|}_{x=x_1+\cdots+x_{r-1}}       \right)  dx_1\cdots dx_{r-2}dx_{r-1}.   \end{split}
    \]
    
    By repeating the above procedure,  for $k-2(r-1)>0$,   it follows that
    \[
    \begin{split}
    &\;\;\;\;\sum_{\substack{k_1+\cdots+k_r=k\\ k_1,\cdots, k_r\geq 1}}f_r(k_1,\cdots,k_{r-1},k_r)   \zeta^{\mathcal{C}}(k_1,\cdots, k_{r-1},1+k_r)\\
    &  =[k-2(r-1)] 
    \sum_{\substack{k_1+\cdots+k_{r-3}+l_{r-2}=k\\k_1,\cdots,k_{r-3}\geq 1,l_{r-2}> 4}}(l_{r-2}-4)\cdots (k_2+\cdots+k_{r-3}+l_{r-2}-2(r-2))\cdot\\
 &\mathop{\int\cdots\int}_{[1,+\infty)^{r-2}}
    \frac{1}{x_1^{k_1}\cdots (x_1+\cdots +x_{r-3})^{k_{r-3}}}\left(\eta\circ \eta\left(\frac{1}{x^{l_{r-2}-3}}\right)\Bigg{|}_{x=x_1+\cdots+x_{r-2}}       \right)  dx_1\cdots dx_{r-3}dx_{r-2}\\ 
    &\;\;\;\;\;\;\;\;\;\;\;\;\; \;\;\;\;\;\;\;\;\;\;\;\;\; \vdots\;\;\;\;\;\;\;\;\;\;\;\;\;\;\;  \;\;\;\;\;\;\;\;\;\;\;\;\;   \vdots  \;\;\;\;\;\;\;\;\;\;\;\;\;   \;\;\;\;\;\;\;\;\;\;\;\;\;\;\;\;\;\;\vdots\\
    &=[k-2(r-1)]\int^{+\infty}_1 \left(  \underbrace{\eta\circ\cdots\circ \eta}_{r-1}\left(\frac{1}{x^{k-2r+1}} \right) \Bigg{|}_{x=x_1}       \right)dx_1 \\
    &=  \underbrace{\eta\circ\cdots\circ \eta}_{r-1}\left(\frac{1}{x^{k-2(r-1)}} \right) \Bigg{|}_{x=1} .       \\
    \end{split}
    \]
    Here the last identity essentially follows from that $\eta$ is a $\mathbb{Q}$-linear transformation on ${\bf V}$.
    As a result, Theorem \ref{sum} is proved.

    Now we give some examples of Theorem   \ref{sum}.
 \begin{ex}
     $(i)$ For $r=2,k>2$, 
     \[
     (k-2)\sum_{\substack{k_1+k_2=k\\k_1,k_2\geq 1}}k_2\zeta^{\mathcal{C}}(k_1,1+k_2)= 1-\frac{1}{2^{k-2}}      ;
     \]
     $(ii)$ For $r=3,k>4$, 
     \[
    (k-4) \sum_{\substack{k_1+k_2+k_3=k\\k_1,k_2,k_3\geq 1}}k_3(k_2+k_3-2)\zeta^{\mathcal{C}}(k_1,k_2, 1+k_3)=\frac{1}{2}-\frac{1}{2^{k-4}}+\frac{1}{2} \frac{1}{3^{k-4}} ;
          \]
      $(iii)$ For $r=4, k>6$,
      \[
      \begin{split}
     &\;\;\;\; (k-6) \sum_{\substack{k_1+k_2+k_3+k_4=k\\k_1,k_2,k_3,k_4\geq 1}} k_4(k_3+k_4-2)(k_2+k_3+k_4-4) \zeta^{\mathcal{C}}(k_1,k_2,k_3,1+k_4)\\
     &=\frac{1}{6}-\frac{1}{2}\cdot \frac{1}{2^{k-6}}+\frac{1}{2}\cdot \frac{1}{3^{k-6}}-\frac{1}{6}\cdot \frac{1}{4^{k-6}}.
     \end{split}    \]
          \end{ex}
          
  \section{Continuous multiple zeta values and multiple polylogarithms}     
         The multiple polylogarithms are defined by
         \[
         \mathrm{Li}_{k_1,\cdots,k_r}(z_1,\cdots, z_r)=\sum_{0<n_1<\cdots<n_r}\frac{z_1^{n_1}\cdots z_r^{n_r}}{n_1^{k_1}\cdots n_r^{k_r}}, \;k_1,\cdots, k_r\geq 1.
         \]  
         For $k_r=1$, $\mathrm{Li}_{k_1,\cdots,k_r}(z_1,\cdots, z_r)$ is convergent for $|z_1|,\cdots, |z_r|<1$. For $k_r\geq 2$, $\mathrm{Li}_{k_1,\cdots,k_r}(z_1,\cdots,z_r)$ is convergent for $|z_1|,\cdots,|z_r|\leq1$. 
         The theory of multiple polylogarithms is related to the $K$-theory of number fields. For $r=1$, they are called polylogarithms.
         
         For a number field $F$, Zagier \cite{zag} conjectured that  Dedekind zeta values $\zeta_F(n)$ can be expressed in terms of polylogarithms at some special algebraic points. For $n=2$, it was proved in \cite{za86}. For $n=3,4$, it was proved in \cite{gon1},\cite{gon2}.
         
         In this section we will show that  a statement which is similar to  Zagier's conjecture holds, for some special continuous multiple zeta values.   
         \begin{lem}\label{fa}
         For $k\geq 1$, one has 
        \[ \frac{1}{x}-\frac{c^k}{x(x+c)^k}=\sum_{0\leq l\leq k-1}\frac{c^l}{(x+c)^{l+1}}.
        \]
         \end{lem}
         \noindent{\bf Proof:} 
         We have 
         \[
        \frac{1}{cx}- \frac{1}{x(x+c)}=\frac{1}{c(x+c)}.
         \]
         Thus for $0\leq l\leq k-1$, one has
         \[\frac{c^l}{x(x+c)^l}-\frac{c^{l+1}}{x(x+c)^{l+1}}=\frac{c^l}{(x+c)^{l+1}}.\]
         It follows that
        \[ \frac{1}{x}-\frac{c^k}{x(x+c)^k}=\sum_{0\leq l\leq k-1}\frac{c^l}{(x+c)^{l+1}}.
        \]
         The lemma is proved. $\hfill\Box$\\            
           
     \begin{Thm}\label{szag}
     Denote by $ \mathfrak{L}_{\mathbb{Q}}$ the $\mathbb{Q}$-algebra generated by 
     the special values of multiple polylogarithms at rational points.
           Then for $r\geq 1$,  $ \zeta^{\mathcal{C}}(\underbrace{1,\cdots,1}_{r},2) \in \mathfrak{L}_{\mathbb{Q}}$.    \end{Thm} 
 \noindent{\bf Proof:}  For $r=1$, the statement is obvious. For $r\geq 2$, by Proposition \ref{algebra}, one has
 \[
 \begin{split}
 &\;\;\;\;\; \zeta^{\mathcal{C}}(\underbrace{1,\cdots,1}_{r},2)\\
 &=\mathop{\int\cdots\int}_{[0,1]^r}\frac{dy_1dy_2\cdots dy_r}{1+y_1+y_1y_2+\cdots+y_1\cdots y_r}\\
 &=\mathop{\int\cdots\int}_{[0,1]^{r-1}} \left( \mathrm{log}\,\left(1+y_1(1+y_2+\cdots+y_2\cdots y_r)\right)\,\bigg{|}^{y_1=1}_{y_1=0}   \right)\frac{ dy_2\cdots dy_r}{1+y_2+\cdots+y_2\cdots y_r}\\
 &=\mathop{\int\cdots\int}_{[0,1]^{r-1}} \frac{\mathrm{log}\,(2+y_2+\cdots+y_2\cdots y_r)}{1+y_2+\cdots+y_2\cdots y_r}dy_2\cdots dy_r.
   \end{split}
 \]
 Since 
 \[
 \begin{split}
 &\;\;\;\;\frac{\mathrm{log}\,(2+y_2+\cdots+y_2\cdots y_r)}{1+y_2+\cdots+y_2\cdots y_r}\\
 &=\frac{\mathrm{log}\,(1+y_2+\cdots+y_2\cdots y_r)}{1+y_2+\cdots+y_2\cdots y_r}+ \frac{\mathrm{log}\,(1+\frac{1}{1+y_2+\cdots+y_2\cdots y_r})}{1+y_2+\cdots+y_2\cdots y_r} \\
 &=\frac{\mathrm{log}\,(1+y_2+\cdots+y_2\cdots y_r)}{1+y_2+\cdots+y_2\cdots y_r}+ \sum_{n\geq 1}\frac{(-1)^{n-1}}{n}\frac{1}{(1+y_2+\cdots+y_2\cdots y_r)^{n+1}},  \end{split}
 \]
 we have
  \[
 \begin{split}
 &\;\;\;\;\; \zeta^{\mathcal{C}}(\underbrace{1,\cdots,1}_{r},2)\\
 &=\mathop{\int\cdots\int}_{[0,1]^{r-2}} \left( \frac{1}{2}\mathrm{log}^2(1+y_2(1+y_3+\cdots+y_3\cdots y_r))\,\bigg{|}^{y_2=1}_{y_2=0}       \right)\frac{dy_3\cdots dy_r}{1+y_3+\cdots+y_3\cdots y_r}\\
 &+\sum_{n\geq 1}\frac{(-1)^{n-1}}{n}\mathop{\int\cdots\int}_{[0,1]^{r-2}} \left( -\frac{1}{n(1+y_2+\cdots+y_2\cdots y_r)^n}\,\bigg{|}^{y_2=1}_{y_2=0}     \right)\frac{dy_3\cdots dy_r}{1+y_3+\cdots+y_3\cdots y_r} \\
 &=\frac{1}{2}\mathop{\int\cdots\int}_{[0,1]^{r-2}}\frac{ \mathrm{log}^2(2+y_3+\cdots+y_3\cdots y_r)}{1+y_3+\cdots+y_3\cdots y_r}dy_3\cdots dy_r\\
 &+\sum_{n\geq 1}\frac{(-1)^{n-1}}{n^2}\mathop{\int\cdots\int}_{[0,1]^{r-2}}\left( 1-\frac{1}{(2+y_3+\cdots+y_3\cdots y_r)^n}  \right) \frac{dy_3\cdots dy_r}{1+y_3+\cdots+y_3\cdots y_r}  \end{split}
 \] 
By Lemma \ref{fa}, we have
\[
\begin{split}
&\;\;\;\;\; \zeta^{\mathcal{C}}(\underbrace{1,\cdots,1}_{r},2)\\
&=\frac{1}{2}\mathop{\int\cdots\int}_{[0,1]^{r-2}}\frac{ \mathrm{log}^2(2+y_3+\cdots+y_3\cdots y_r)}{1+y_3+\cdots+y_3\cdots y_r}dy_3\cdots dy_r\\
&\;\;\;\;\;\;\;+\sum_{1\leq n_1\leq n}\frac{(-1)^{n-1}}{n^2}\mathop{\int\cdots\int}_{[0,1]^{r-2}} \frac{dy_3\cdots dy_r}{(2+y_3+\cdots+y_3\cdots y_r)^{n_1}}\\
&=\frac{1}{2}\mathop{\int\cdots\int}_{[0,1]^{r-2}}\frac{\mathrm{log}^2(2+y_3+\cdots+y_3\cdots y_r)}{1+y_3+\cdots+y_3\cdots y_r}dy_3\cdots dy_r\\
&+\sum_{n\geq 1}\frac{(-1)^{n-1}}{n^2}\mathop{\int\cdots\int}_{[0,1]^{r-2}} \frac{dy_3\cdots dy_r}{2+y_3+\cdots+y_3\cdots y_r}\\
&+\sum_{1\leq n_1<n}\frac{(-1)^{n-1}}{n^2}\mathop{\int\cdots\int}_{[0,1]^{r-2}} \frac{dy_3\cdots dy_r}{(2+y_3+\cdots+y_3\cdots y_r)^{n_1+1}}.\\
\end{split} \tag{1}
\]
From the above analysis, we have 
\[
\zeta(1,1,2)=\frac{\mathrm{log}^22}{2}+\sum_{n\geq 1}\frac{(-1)^{n-1}}{n^2}(1-\frac{1}{2^n})=\frac{(\mathrm{Li}_1(-1))^2}{2}-\mathrm{Li}_2(-1)+\mathrm{Li}_2\left(-\frac{1}{2}\right)
\]
and
\[
\begin{split}
&\;\;\;\;\zeta(1,1,1,2)\\
&=\frac{1}{2}\int_0^1\frac{\mathrm{log}^2(2+y)}{1+y}dy+\sum_{n\geq 1}\frac{(-1)^{n-1}}{n^2}\left(\int^1_0\frac{dy}{2+y}+\int^1_0\frac{dy}{(2+y)^{n+1}}\right)\\
&=\frac{1}{2}\int_0^1\frac{\mathrm{log}^2(2+y)}{1+y}dy+\sum_{n\geq 1}\frac{(-1)^{n-1}}{n^2}\left[\mathrm{log}\,\frac{3}{2}+\frac{1}{n}\left(\frac{1}{2^n}-\frac{1}{3^n}\right)\right]\\
&=\frac{1}{2}\int_0^1\frac{1}{1+y}\left[ \mathrm{log}^2(1+y)+ 2\mathrm{log}(1+y)\sum_{n\geq1}\frac{(-1)^{n-1}}{n}\frac{1}{(1+y)^n}+\sum_{n_1,n_2\geq 1}\frac{(-1)^{n_1+n_2}}{n_1n_2}\frac{1}{(1+y)^{n_1+n_2}}      \right]          dy\\
&\;\;\;\;-\mathrm{Li}_2(-1)\mathrm{Li}_1\left(\frac{1}{3}\right)-\mathrm{Li}_3\left(-\frac{1}{2} \right)+\mathrm{Li}_3\left(-\frac{1}{3} \right)\\
&=\frac{1}{6}\mathrm{log}^32+\sum_{n\geq 1}\frac{(-1)^{n-1}}{n}\int^1_0\frac{\mathrm{log}(1+y)}{(1+y)^{n+1}}dy+\frac{1}{2}\sum_{n_1,n_2\geq 1}\frac{(-1)^{n_1+n_2}}{n_1n_2(n_1+n_2)}\left(1-\frac{1}{2^{n_1+n_2}}\right)\\
&\;\;\;\;-\mathrm{Li}_2(-1)\mathrm{Li}_1\left(\frac{1}{3}\right)-\mathrm{Li}_3\left(-\frac{1}{2} \right)+\mathrm{Li}_3\left(-\frac{1}{3} \right).\\
\end{split}
\]
By integration by parts, one can check that
\[
\int^1_0\frac{\mathrm{log}(1+y)}{(1+y)^{n+1}}dy=- \frac{1}{n}\frac{\mathrm{log}\,2}{2^n}+\frac{1}{n^2}\left(1-\frac{1}{2^n}\right).
\]
Thus we have 
\[
\begin{split}
& \;\;\;\;\zeta(1,1,1,2)  \\
&=\frac{1}{6}\mathrm{log}^32+\sum_{n\geq 1}\frac{(-1)^{n-1}}{n}\left[ - \frac{1}{n}\frac{\mathrm{log}\,2}{2^n}+\frac{1}{n^2}\left(1-\frac{1}{2^n}\right)       \right]      +\frac{1}{2}\sum_{n_1,n_2\geq 1}\frac{(-1)^{n_1+n_2}}{n_1n_2(n_1+n_2)}\left(1-\frac{1}{2^{n_1+n_2}}\right)\\
&\;\;\;\;-\mathrm{Li}_2(-1)\mathrm{Li}_1\left(\frac{1}{3}\right)-\mathrm{Li}_3\left(-\frac{1}{2} \right)+\mathrm{Li}_3\left(-\frac{1}{3} \right)\\
&=\frac{1}{6}\mathrm{log}^32+\mathrm{Li}_2\left(-\frac{1}{2}\right)\mathrm{log}\,2-\mathrm{Li}_3(-1)+\mathrm{Li}_3\left(-\frac{1}{2}\right)+\sum_{n_1,n_2\geq 1}\frac{(-1)^{n_1+n_2}}{n_1(n_1+n_2)^2}\left(1-\frac{1}{2^{n_1+n_2}}\right) \\
&\;\;\;\;-\mathrm{Li}_2(-1)\mathrm{Li}_1\left(\frac{1}{3}\right)-\mathrm{Li}_3\left(-\frac{1}{2} \right)+\mathrm{Li}_3\left(-\frac{1}{3} \right)\\
&=-\frac{(\mathrm{Li}_{1}(-1))^3}{6}-\mathrm{Li}_2\left(-\frac{1}{2}\right)\mathrm{Li}_1(-1)-\mathrm{Li}_3(-1)+\mathrm{Li}_3\left(-\frac{1}{2}\right)+    \mathrm{Li}_{1,2}(-1)-\mathrm{Li}_{1,2}\left(-\frac{1}{2}\right)\\
&\;\;\;\;-\mathrm{Li}_2(-1)\mathrm{Li}_1\left(\frac{1}{3}\right)-\mathrm{Li}_3\left(-\frac{1}{2} \right)+\mathrm{Li}_3\left(-\frac{1}{3} \right).\\
\end{split}
\]
For $r>3$, one can use the  formula $(1)$, Lemma \ref{fa} and the trick
\[
\begin{split}
&\;\;\;\;\mathrm{log}(k+1+t_1+\cdots+t_1\cdots t_r)\\
&=\mathrm{log}(k+t_1+\cdots+t_1\cdots t_r)+\mathrm{log}\left(  1+\frac{1}{k+t_1+\cdots+t_1\cdots t_r}\right)\\
&=\mathrm{log}(k+t_1+\cdots+t_1\cdots t_r)+\sum_{n\geq 1}\frac{(-1)^{n-1}}{n}\frac{1}{(k+t_1+\cdots+t_1\cdots t_r)^n}.
\end{split}
\]
repeatedly to deduce that
$$ \zeta^{\mathcal{C}}(\underbrace{1,\cdots,1}_{r},2) \in \mathfrak{L}_{\mathbb{Q}}.$$

In a word, from the above analysis, we give an explicit procedure to calculate $$ \zeta^{\mathcal{C}}(\underbrace{1,\cdots,1}_{r},2)$$ in terms of elements in $\mathfrak{L}_{\mathbb{Q}}$.
$\hfill\Box$\\     

\begin{rem}
In general cases, by changing of variables
\[
t_i=\frac{1}{x_i}, i=1,\cdots,x_r,
\]
one has  \[
\begin{split}
&\;\;\;\zeta^{\mathcal{C}}(k_1,k_2,\cdots, k_r)\\
&=\mathop{\int\cdots \int}_{[1,+\infty)^r}\frac{dx_1dx_2\cdots dx_r}{x_1^{k_1}(x_1+x_2)^{k_2}\cdots (x_1+\cdots+x_r)^{k_r}}\\
&=\mathop{\int}_{[0,1]^r} \frac{1}{(\frac{1}{t_1})^{k_1} (\frac{1}{t_1}+\frac{1}{t_2})^{k_2}\cdots (  \frac{1}{t_1}+\cdots+\frac{1}{t_r}    )^{k_r} }\frac{dt_1dt_2\cdots dt_r}{t_1^2\;t_2^2\;\cdots \;t_r^2} . \\
\end{split}
\]
Thus all the continuous multiple zeta values are periods in the sense of Kontsevich-Zagier \cite{KZ}. For now, we only know that $\mathrm{log}\,2$ and $\mathrm{log}\,3$(which are elements in the algebra of continuous multiple zeta values) are transcendental. Theorem \ref{szag} shows that the algebra of continuous multiple zeta values is closely related to the special values of multiple polylogarithms.
\end{rem}       
                       
          \section{Further remarks}
            In this section we discuss some related topics which are still unclear for now.
                      \subsection{Depth defect phenomena for continuous multiple zeta values}
          From Section \ref{de} we have that
          \[\mathrm{dim}_{\mathbb{Q}} \mathfrak{D}_r\mathcal{Z}^{\mathcal{C}}= 2^{r-1}, r=1,2,\]
               \[\mathrm{dim}_{\mathbb{Q}} \mathfrak{D}_3\mathcal{Z}^{\mathcal{C}}\leq 4, \;\mathrm{dim}_{\mathbb{Q}} \mathfrak{D}_r\mathcal{Z}^{\mathcal{C}}< 2^{r-1}, r\geq 4.\]
               Thus there are depth defect phenomena for continuous multiple zeta values of depth $\geq 4$. Can we give a sharper bound for $\mathrm{dim}_{\mathbb{Q}} \mathfrak{D}_r\mathcal{Z}^{\mathcal{C}}$ in case $r\geq 4$?

          \subsection{Continuous multiple zeta values and cyclotomic multiple zeta values}
          For $N\geq 1$, cyclotomic multiple zeta values of level $N$ are defined by 
          \[
          \zeta\binom{k_1,\cdots,k_r}{\epsilon_1,\cdots,\epsilon_r}=\sum_{0<n_1<\cdots<n_r}\frac{\epsilon_1^{n_1}\cdots \epsilon_r^{n_r}}{n_1^{k_1}\cdots n_r^{k_r}},
          \]
          where $k_1,\cdots,k_r\geq1, (k_r,\epsilon_r)\neq (1,1),\epsilon_1^N=\cdots =\epsilon_r^N=1$.
          The set of $\mathbb{Q}$-linear combinations of cyclotomic multiple zeta values of level $N$ is also a $\mathbb{Q}$-algebra.
          
          It is clear that $\mathrm{log}\;2$ is a cyclotomic multiple zeta value of level $2$. 
          More generally, by the decomposition
          \[
          x^n-1=\prod_{\epsilon\in \mu_N}(x-\epsilon)     ,\; \mu_N=\{\epsilon\,|\,\epsilon^N=1\},
          \] 
          one can check that $\mathrm{log}\;N$  is a cyclotomic multiple zeta value of level $N$ for $N\geq 2$. 
          
          Since $\mathrm{log}\;2, \mathrm{log}\;3$ are both continuous multiple zeta values, the algebra of continuous multiple zeta values and the algebra of cyclotomic multiple zeta values of level $N$ have some common elements for some $N$. It is very interesting to give a detailed
           analysis about the intersection of the above two sets.
                  
\section*{Acknowledgements}
The author wants to thank the anonymous referee for his/her helpful comments to improve this paper. The author is supported by the National Natural Science Foundation of China (Grant No. 12201642).


\begin{thebibliography}{1}


\bibitem{bk} 
   D. Broadhurst, D. Kreimer,
   \emph{Association of multiple zeta values with positive knots via Feynman diagrams up to $9$ loops},
    Phys. Lett. B 393 (1997), no. 3-4, 303-412.

\bibitem{chen}
	K.  Chen,
	\emph{Iterated path integrals},
	Bull. Amer. Math. Soc.,(1977), 83, 831-879.

\bibitem{gkz}
    H. Gangl, M. Kaneko, D. Zagier,
    \emph{Double zeta values and modular forms},
    Automorphic forms and zeta functions,
    In: Proceedings of the conference in memory of Tsuneo Arakawa,
    World Scientific (2006),
    71-106.
\bibitem{gon1}
A. Goncharov, Geometry of configurations, polylogarithms, and motivic cohomology,
Adv. Math. 114 (1995), no. 2, 197-318.

\bibitem{gon2}
A. Goncharov, D. Rudenko,  \emph{Motivic correlators, cluster varieties, and Zagier’s conjecture on $\zeta_F(4)$}, arXiv: 1803.08585.

\bibitem{gon}
    A. Goncharov, \emph{The dihedral Lie algebras and Galois symmetries of $\pi_1^{(l)}(\mathbb{P}^1-\{0,\infty\}\cup \mu_{N})$},
    Duke Math. J., 110 (3)  (2001), 397-487.


\bibitem{ikz}

K. Ihara, M. Kaneko, D. Zagier,
 \emph{Derivation and double shuffle relations for multiple zeta values},
Compositio Math. 142 (2006), 307-338.

\bibitem{KZ}
M. Kontsevich, D. Zagier, \emph{Periods}, Mathematics Unlimited 2001 and Beyond, Springer, Berlin, (2001), 771-808.

\bibitem{yam}
S. Yamamoto, \emph{Multiple zeta-star values and multiple integrals}, arXiv: 1405.6499.

\bibitem{za86}
D. Zagier, \emph{Hyperbolic manifolds and special values of Dedekind zeta functions}, Inventiones Math., 83 (1986), 285-301. 
\bibitem{zag}
D. Zagier, \emph{Polylogarithms, zeta-functions, and algebraic K-theory of fields}, Progress
Math Vol 89. Birkhauser, Boston, MA, (1991), 392-430.



\bibitem{zhao}
J. Zhao, \emph{Analytic continuation of multiple zeta functions }, Proceedings of the American Mathematical Society, 128 (5) (1999), 1275-1283.


\end{thebibliography}
\end{document}